\documentclass[12pt,a4paper,twoside]{amsart}
\usepackage{amssymb}
\usepackage{amsmath}
\usepackage{amsfonts}
\usepackage{bm}
\usepackage{graphicx}
\usepackage{mathrsfs}
\usepackage{fancyhdr}
\usepackage{subfigure}  
\usepackage[mathscr]{eucal}
\usepackage{amsxtra}
\usepackage{latexsym}
\usepackage{multicol}
\usepackage{fancybox}
\usepackage{helvet}
\usepackage[T1]{fontenc}
\usepackage{wasysym}
\usepackage{colortbl}
\usepackage{url}

  \usepackage{latexsym}
  \usepackage{amsthm}

\newcommand{\field}[1]{\mathbb{#1}}
\newcommand{\C}{\field{C}}
\newcommand{\R}{\field{R}}  
\newcommand{\Z}{\field{Z}} 

\newcommand{\N}{\field{N}}

\newcommand{\HQ}{\field{H}}
\newcommand{\vect}[1]{\mbox{\boldmath $#1$}}

 \newtheorem{thm}{Theorem}[section]
 \newtheorem{lem}[thm]{Lemma}
 \newtheorem{defn}[thm]{Definition}

  \newcommand{\qi}{\ensuremath{\mbox{\boldmath $i$}}}
  \newcommand{\sqi}{\ensuremath{\mbox{\boldmath $\scriptstyle i$}}}
  \newcommand{\qj}{\ensuremath{\mbox{\boldmath $j$}}}
  \newcommand{\sqj}{\ensuremath{\mbox{\boldmath $\scriptstyle j$}}}
  \newcommand{\qk}{\ensuremath{\mbox{\boldmath $k$}}}

  \newcommand{\be}{\begin{equation}}
  \newcommand{\ee}{\end{equation}} 

  \newcommand{\bfr}{\begin{frame}}
  \newcommand{\efr}{\end{frame}}

\usepackage{calligra}
\DeclareMathAlphabet{\mathcalligra}{T1}{calligra}{m}{n}
\DeclareFontShape{T1}{calligra}{m}{n}{<->s*[2.2]callig15}{}

\begin{document}
         
\title{Quaternionic Fourier-Mellin Transform}

\author{Eckhard Hitzer}
\address{Department of Applied Physics, University of Fukui,\\ 
Japan}
\email{hitzer@mech.u-fukui.ac.jp}

\maketitle

\lhead{{\scriptsize Proceedings of the \\
19th ICFIDCAA Hiroshima 2011 \\ 
Tohoku UP, 2012, pp.~123--133}}
\chead{}
\rhead{}
\renewcommand{\headrulewidth}{0pt}
\thispagestyle{fancy}

\begin{abstract}
In this contribution we generalize the classical Fourier Mellin transform \cite{DG:RobbEffFMT}, which transforms functions $f$ representing, e.g., a gray level image defined over a compact set of $\mathbb{R}^2$. 
The quaternionic Fourier Mellin transform (QFMT) applies to functions
$f: \mathbb{R}^2 \rightarrow \mathbb{H}$, for which $|f|$ is summable over $\mathbb{R}_+^* \times \mathbb{S}^1$ under the measure $d\theta \frac{dr}{r}$. $\mathbb{R}_+^*$ is the multiplicative group of positive and non-zero real numbers.
We investigate the properties of the QFMT similar to the investigation of the quaternionic Fourier Transform (QFT) in \cite{EH:QFTgen,EH:DirUP}. 
\end{abstract}

\section{Quaternions}


Gauss\footnote{Remark:
By reading this paper you agree to the terms of the \textit{Creative Peace License} of page \pageref{pg:CPL}.}, Rodrigues and Hamilton introduced the {4D quaternion algebra} $\HQ$ over $\R$ with {three imaginary units}:
\be
 \qi \qj = -\qj \qi = \qk, \,\,
 \qj \qk = -\qk \qj = \qi, \,\,
 \qk \qi = -\qi \qk = \qj, \,\, 
 \qi^2=\qj^2=\qk^2=\qi \qj \qk = -1.
\label{eq:quat}
\end{equation}
Every quaternion 
\be
  q=q_r + q_i \qi + q_j \qj + q_k \qk \in \HQ, \quad 
  q_r,q_i, q_j, q_k \in \R
  \label{eq:aquat}
\end{equation}
has {\textit{quaternion conjugate}} (reversion in $Cl_{3,0}^+$)
\be
  \tilde{q} = q_r - q_i \qi - q_j \qj - q_k \qk,
\end{equation}
This leads to {\textit{norm}} of $q\in\HQ$, and an {inverse} of every non-zero $q\in\HQ$
\be
  | q | = \sqrt{q\tilde{q}} = \sqrt{q_r^2+q_i^2+q_j^2+q_k^2},
  \quad
  | p q | = | p || q |,
  \quad
  q^{-1} = \frac{\tilde{q}}{|q|^2} = \frac{\tilde{q}}{q\tilde{q}} . 
\end{equation}
The scalar part of quaternions is symmetric
\be
  Sc(q) = q_r = \frac{1}{2}(q+\tilde{q}), 
  \quad 
  Sc(pq) = Sc(qp) . 
\ee
The \textit{inner product} of quaternions defines \textit{orthogonality} 
\be 
  Sc(p\widetilde{q}) 
  = p_rq_r + p_iq_i + p_jq_j + p_kq_k \in \R . 
  \label{eq:4Dinnp}
\ee

\subsection{The (2D) orthogonal planes split (OPS) of quaternions}

We consider an arbitrary pair of pure quaternions $f,g$, $f^2=g^2=-1$. 
The {orthogonal 2D planes split (OPS)} is then defined with respect to a pair of pure unit quaternions $f, g$ as 
\be 
q_{\pm} = \frac{1}{2}(q \pm f q g).
\ee 

We thus observe, that 
\be 
  f q g = q_+ - q_-,
\ee 
i.e. under the map $f( )g$ the {$q_+$ part is invariant}, but the {$q_-$ part changes sign}.
Both parts are two-dimensional, and span {two completely orthogonal planes}. 
For $f\neq \pm g$ the {$q_+$ plane} is spanned by the orthogonal quaternions $\{f-g, 1+fg\}$,
and the {$q_-$ plane} is spanned by $\{f+g, 1-fg\}$. 
Vice versa, in general any two 2D orthogonal planes in $\HQ$ determine a corresponding pair $f,g$.

\begin{lem}[Orthogonality of two OPS planes]
\label{lm:OPSortho}
  Given any two quaternions $q,p$ and applying the OPS with respect to two linearly independent pure unit quaternions $f,g$ we get zero for the
  scalar part of the mixed products
  \be
     Sc(p_+\widetilde{q}_-) = 0, \quad Sc(p_-\widetilde{q}_+) = 0 
     \quad \Longrightarrow \quad
     |q|^2 = |q_+|^2 + |q_-|^2 .
  \ee
\end{lem}

\subsection{Geometric interpretation of map $f()g$}

\begin{figure}
\begin{center}
  \includegraphics[height=4.0cm]{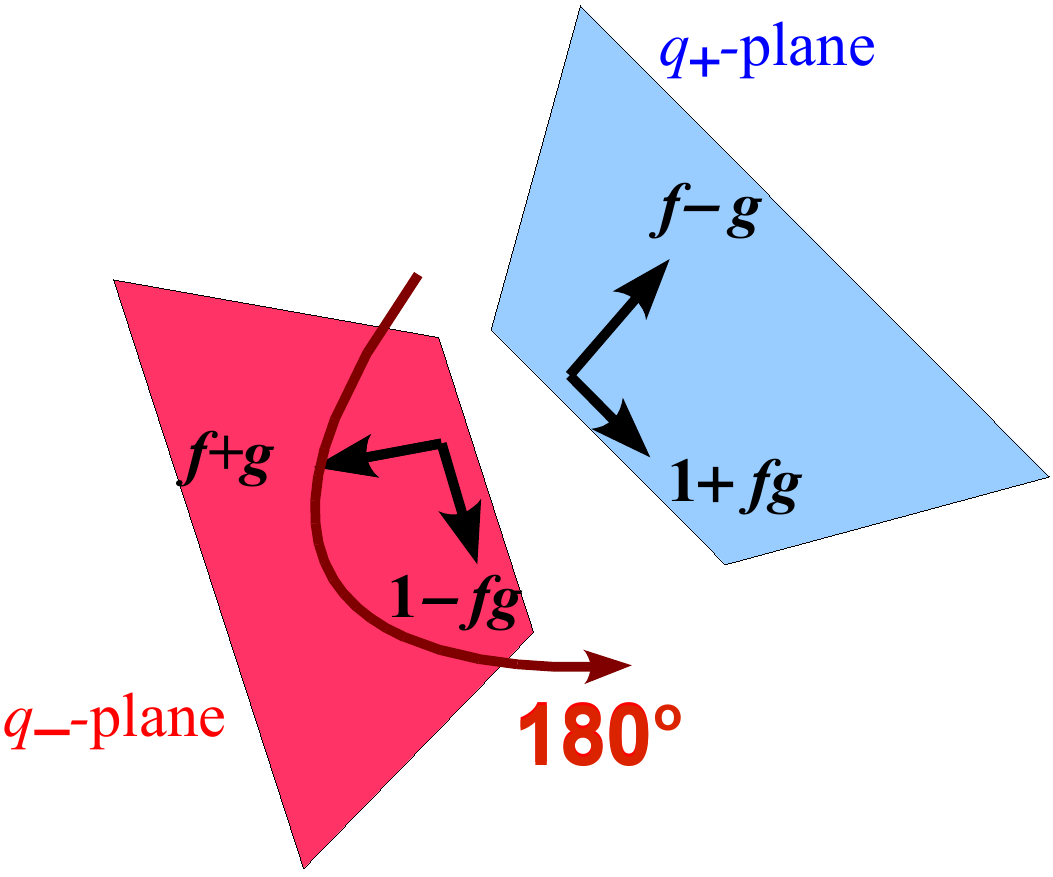}
  \caption{The map $f( )g$ {rotates} the $q_-$ plane by $180^{\circ}$ around the $q_+$ axis plane. (Basis for $g\neq \pm f$.) \label{fg:HalfTurn}}
\end{center}
\end{figure}

The map $f( )g$ {rotates} the $q_-$ plane by $180^{\circ}$ around the $q_+$ axis plane, see Fig.\footnote{Color versions of the figures (showing maxima and minima in red and blue colors, respectively), can be found in the electronic copy of this contribution to appear on \url{http://sinai.mech.fukui-u.ac.jp/gcj/pubs.html} .} \ref{fg:HalfTurn}. This interpretation of the map $f( )g$ is in perfect agreement with Coxeter's notion of {half-turn} in \cite{HSMC:QuatRef}. 

We obtain the following important identities:
\be 
  e^{\alpha f} q_{\pm} e^{\beta g} 
  = q_{\pm} e^{(\beta\mp\alpha) g}
  = e^{(\alpha\mp\beta) f}q_{\pm},
  \qquad 
  q\in \HQ, \,\, \alpha, \beta \in \R.
  \label{eq:eqpmeres}
\ee 

\begin{figure}
\begin{center}
  \includegraphics[height=5.5cm]{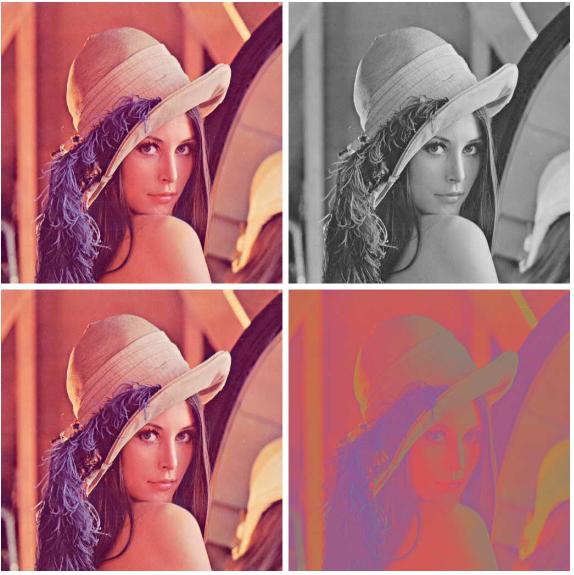}
  \caption{
  {Simplex/perplex split} of \cite{ES:HFTColIm} with gray line $f=g=(\qi+\qj+\qk)/\sqrt{3}$. 
  \textit{Top left}: Original. 
  \textit{Top right}: $q_-$-part (luminance). 
  \textit{Bottom right}: $q_+$-part ($q_i \leftrightarrow q_j$) (chrominance). 
  \textit{Bottom left}: Sum. \label{fg:lumchrom}}
\end{center}
\end{figure}

For $g\neq \pm f$ the set $\{f-g, 1+fg, f+g, 1-fg\}$ forms an {orthogonal basis of $\HQ$} interpreted as $\R^4$. 
We can therefore use the following {representation} for every $q \in \HQ$ by means of four real coefficients $q_1, q_2, q_3, q_4 \in \R$
\begin{gather} 
  q = q_1 (1+fg) + q_2 (f-g) + q_3 (1-fg) + q_4 (f+g),
  \\
  q_1 = Sc(q(1+fg)^{-1}), \quad 
  q_2 = Sc(q(f-g)^{-1}), \nonumber \\
  q_3 = Sc(q(1-fg)^{-1}), \quad
  q_4 = Sc(q(f+g)^{-1}).
  \nonumber 
\end{gather}
In the case of $f=\qi, g=\qj$ we obtain the coefficients
\be 
  q_1 = \frac{1}{2}(q_r+q_k), \quad 
  q_2 = \frac{1}{2}(q_i-q_j), \quad
  q_3 = \frac{1}{2}(q_r-q_k), \quad
  q_4 = \frac{1}{2}(q_i+q_j).
\ee 
The OPS with respect to a {single pure unit quaterion}, 
e.g., $f=g=\qi$ gives
\be 
  q_{\pm} = \frac{1}{2}(q\pm \qi q \qi), \quad 
  q_+ = q_j \qj + q_k \qk = (q_j + q_k \qi)\qj, \quad 
  q_- = q_r + q_i \qi ,
  \label{eq:opsiqi}
\ee  
where the $q_+$ plane is 2D and manifestly orthogonal to the 2D $q_-$ plane. 
The above corresponds to the {simplex/perplex split} of \cite{TAE:QFT}, see an application in Fig. \ref{fg:lumchrom} from \cite{ES:HFTColIm}.

\section{The Quaternionic Fourier Mellin transformations (QFMT)}

\subsection{Robert Hjalmar Mellin (1854--1933)}

\begin{figure}
\begin{center}
  \includegraphics[height=3.0cm]{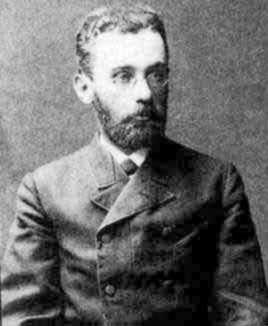}
  \caption{Robert Hjalmar Mellin (1854--1933). Image: Wikipedia. \label{fg:RHMellin}}
  \end{center}
\end{figure}

Robert Hjalmar Mellin (1854--1933) \cite{StA:Mellinbio}, Fig. \ref{fg:RHMellin}, was a Finnish mathematician, a student of G. Mittag-Leffler and K. Weierstrass. He became the director of the Polytechnic Institute in Helsinki, and in 1908 first professor of mathematics at Technical University of Finland. He was a fervent fennoman with fiery temperament, and co-founder of the Finnish Academy of Sciences. He became known for the \textit{Mellin transform} with major applications to the evaluation of integrals, see \cite{PBM:EvIntMellT}, which lists {1624 references}. During his last 10 years he tried to {refute Einstein's theory} of relativity as logically untenable. 

\begin{defn}[Definition: Classical Fourier Mellin transform (FMT)]
\begin{equation}
  \forall(v,k)\in \mathbb{R} \times \mathbb{Z}, \quad
  \mathcal{M}\{h\}(v,k) = \frac{1}{2\pi}\int_0^{\infty}\int_0^{2\pi}
  h(r,\theta) r^{-iv}e^{-ik\theta}d\theta \frac{dr}{r},
\end{equation}
where $h: \mathbb{R}^2 \rightarrow \R$ denotes a function representing, e.g., a gray level image defined over a compact set of $\mathbb{R}^2$. 
\end{defn}
Well known applications are to shape recognition (independent of rotation and scale), image registration, and similarity. 

\subsection{Inner product, symmetric part, norm of quaternion-valued functions}

The quaternion $\HQ$-valued inner product for quaternion-valued functions
$h,m: \R^2 \rightarrow \HQ$ is given by
\begin{equation}
  (h,m) = \int_{\R^2} h(\vect{x})\,\widetilde{m}(\vect{x})\,d^2\vect{x} \,, 
  \quad \quad
  \text{ with } \quad d^2\vect{x} = dx dy,
\label{eq:intip}
\end{equation}
It has a symmetric real {scalar part}
\begin{equation}
  \langle h,m\rangle = \frac{1}{2}[(h,m)+(m,h)] 
  = \int_{\R^2} Sc\left( h(\vect{x})\,\widetilde{m}(\vect{x}) \right) d^2\vect{x} \, ,
\label{eq:intsp}
\end{equation}  
which allows to define a $L^2(\R^2;\HQ)$-norm
\begin{equation}
  \| h \|^2 = (h,h)
          = \langle h,h\rangle
          = \int_{\R^2} |h(\vect{x})|^2 \,d^2\vect{x}
  \qquad 
  \Longrightarrow \,\,\, \|h\|^2 = \|h_+\|^2 + \|h_-\|^2 \,\,.
\end{equation}
A quaternion module can be defined as $L^2(\R^2;\HQ)$ by
\begin{equation}
  L^2(\R^2;\HQ ) = \{h | h:\R^2 \rightarrow \HQ, \|h\| < \infty\}.
\end{equation}

We now define the generalization of the FMT to quaternionic signals. 

\begin{defn}[Quaternionic Fourier Mellin transform (QFMT)]
Let $f,g \in \HQ: f^2=g^2=-1$ be any pair of pure unit quaternions. 
The quaternionic Fourier Mellin transform (QFMT) is given by
\begin{equation}
  \forall(v,k)\in \mathbb{R} \times \mathbb{Z},\quad
  \hat{h}(v,k)=
  \mathcal{M}\{h\}(v,k) 
  = \frac{1}{2\pi}\int_0^{\infty}\int_0^{2\pi}
  r^{-f v} h(r,\theta) e^{-g k\theta}d\theta \frac{dr}{r},
\end{equation}
where $h: \mathbb{R}^2 \rightarrow \mathbb{H}$ denotes a function from $\mathbb{R}^2$ into the algebra of quaternions $\mathbb{H}$, 
such that $|h|$ is summable over $\mathbb{R}_+^* \times \mathbb{S}^1$ under the measure $d\theta \frac{dr}{r}$. $\mathbb{R}_+^*$ is the multiplicative group of positive and non-zero real numbers.
\end{defn}
For $f=\qi$, $g=\qj$ we have the special case
\begin{equation}
  \forall(v,k)\in \mathbb{Z}\times \mathbb{R},\quad
  \hat{h}(v,k)=
  \mathcal{M}\{h\}(v,k) = \frac{1}{2\pi}\int_0^{\infty}\int_0^{2\pi}
  r^{-\sqi v} h(r,\theta) e^{-\sqj k\theta}d\theta \frac{dr}{r},
\end{equation}
Note, that the $\pm$ split and the QFMT {commute}: 
$$
  \mathcal{M}\{h_{\pm}\} = \mathcal{M}\{h\}_{\pm}.
$$

\begin{thm}[Inverse QFMT]
The QFMT can be inverted by
\be 
  h(r,\theta) = 
  \mathcal{M}^{-1}\{h\}(r,\theta)
  = \frac{1}{2\pi}\int_{-\infty}^{\infty}\sum_{k\in \Z}
  r^{f v} \,\hat{h}(v,k) \,e^{g k\theta} dv.
\ee
\end{thm}
The proof uses
\be
  \frac{1}{2\pi}\sum_{k\in \Z}e^{g k(\theta-\theta')} 
  = \delta(\theta-\theta'),
  \qquad
  r^{f v} = e^{f v \ln r}, 
  \qquad
  \frac{1}{2\pi}\int_0^{2\pi}e^{f v (\ln (r) - s)} dv 
  = \delta(\ln (r) - s) .
\ee

We now investigate the basic properties of the QFMT. 
First, left linearity:
For $\alpha, \beta \in \{q \mid q = q_r + q_f f, \, q_r, q_f \in \R\}$,
\be 
  m(r,\theta) = \alpha h_1(r,\theta) + \beta h_2(r,\theta)
  \,\,\, \Longrightarrow \,\,\, 
  \hat{m}(v,k) = \alpha \hat{h}_1(v,k) + \beta \hat{h}_2(v,k) . 
\ee 

Second, right linearity:
For $\alpha', \beta' \in \{q \mid q = q_r + q_g g, \, q_r, q_g \in \R\}$,
\be 
  m(r,\theta) = h_1(r,\theta)\alpha' + h_2(r,\theta)\beta'
  \,\,\, \Longrightarrow \,\,\, 
  \hat{m}(v,k) = \hat{h}_1(v,k)\alpha' + \hat{h}_2(v,k)\beta' . 
\ee 

The linearity of the QFMT leads to
\be
  \mathcal{M}\{h\}(v,k)
  = \mathcal{M}\{h_- + h_+\}(v,k)
  = \mathcal{M}\{h_-\}(v,k)
    + \mathcal{M}\{ h_+\}(v,k),
\ee
which gives rise to the following thoerem.

\begin{thm}[Quasi-complex FMT like forms for QFMT of $h_{\pm}$]
\label{th:fpmtrafo}
The QFMT of $h_{\pm}$ parts of $h \in L^2(\R^2,\HQ)$ 
have simple {quasi-complex forms}
\begin{equation}
 \mathcal{M}\{h_{\pm}\} 
  \stackrel{}{=} \frac{1}{2\pi}\int_0^{\infty}\int_0^{2\pi}
    h_{\pm} r^{\pm g v} e^{-g k \theta} d\theta \frac{dr}{r}
  \stackrel{}{=} \frac{1}{2\pi}\int_0^{\infty}\int_0^{2\pi}
    r^{-f v} e^{\pm f k \theta} h_{\pm} d\theta \frac{dr}{r} \,\, .
\end{equation}
\end{thm}
Theorem \ref{th:fpmtrafo} allows to use discrete and fast software to compute the QFMT based on a pair of complex FMT transformations. 

For the two split parts of the QFMT, we have the following lemma. 
\begin{lem}[Modulus identities]
Due to $|q|^2 = |q_-|^2 + |q_+|^2$ we get for $f:\R^2\rightarrow \HQ$ the following identities
\begin{align}
  |h(r,\theta)|^2 
  &= |h_-(r,\theta)|^2 + |h_+(r,\theta)|^2,
  \nonumber \\
  |\mathcal{M}\{h\}(v,k)|^2 
  &= |\mathcal{M}\{h_-\}(v,k)|^2 + |\mathcal{M}\{h_+\}(v,k)|^2.
\end{align}
\end{lem}

Further properties are \textit{scaling} and \textit{rotation}:
For 
$ m(r,\theta) = h(ar,\theta+\phi)$, 
$\,a > 0, \,\, 0\leq \phi \leq 2\pi $,
\be
  \widehat{m}(v,k) = a^{f v} \hat{h}(v,k) e^{g k \phi} .
\ee 

Moreover, we have the following magnitude identity:
\be 
  |\widehat{m}(v,k)| = |\hat{h}(v,k)| ,
  \label{eq:magid}
\ee 
i.e. the magnitude of the QFMT of a scaled and rotated quaternion signal 
$m(r,\theta) = h(ar,\theta+\phi)$ is identical to the magnitude of the QFMT of $h$.
Equation (\ref{eq:magid}) forms the basis for applications to {rotation and scale invariant shape recognition} and image registration. This may now be extended to {color images}, since quaternions can encode colors RGB in their $\qi, \qj, \qk$ components. 

The reflection at the unit circle ($r \rightarrow \frac{1}{r}$) leads to
\be 
  m(r,\theta) = h(\frac{1}{r}, \theta) 
  \qquad \Longrightarrow \qquad
  \widehat{m}(v,k) = \hat{h}(-v,k) .
\ee

Reversing the sense of sense of rotation ($\theta \rightarrow -\theta$) yields
\be 
  m(r,\theta) = h(r, -\theta) 
  \qquad \Longrightarrow \qquad
  \widehat{m}(v,k) = \hat{h}(v,-k) .
\ee

Regarding radial and rotary modulation we assume
\be
  m(r,\theta) = r^{f v_0} \,h(r,\theta) \,e^{g k_0\theta}, 
  \qquad
  v_0 \in \R, \, k_0 \in \Z . 
\ee
Then we get
\be
  \widehat{m}(v,k) = \hat{h}(v-v_0, k-k_0) .
\ee

\subsection{QFMT derivatives and power scaling}

We note for the logarithmic derivative that $\frac{d}{d \ln r} = r \frac{d}{dr} = r \partial_r$, 
\be 
  \mathcal{M}\{(r\partial_r)^n h\}(v,k) = (f v)^n \hat{h}(v,k), 
  \qquad n \in \N .
\ee
Applying the angular derivative with respect to $\theta$ we obtain
\be 
  \mathcal{M}\{\partial_{\theta}^n h\}(v,k) = \hat{h}(v,k) (g k)^n , 
  \qquad n \in \N .
\ee
Finaly, power scaling with $\ln r$ and $\theta$ leads to
\be 
  \mathcal{M}\{ (\ln r)^m \theta^n h\}(v,k)
  = f^m \,\partial_v^m \partial_k^n \hat{h}(v,k) \,g^n ,
  \qquad
  m,n \in \N .
\ee

\subsection{QFMT Plancherel and Parseval theorems}

For the QFMT we have the following two theorems. 

\begin{thm}[QFMT Plancherel theorem]
The scalar part of the inner product of two functions $h,m : \R^2 \rightarrow \HQ$
is
\be
  \langle h,m \rangle = \langle \hat{h},\widehat{m} \rangle .
\ee 
\end{thm}
\begin{thm}[QFMT Parseval theorem]
Let $h : \R^2 \rightarrow \HQ$. Then
\be
  \|h\| = \|\hat{h}\| ,
  \qquad
  \|h\|^2 = \|\hat{h}\|^2 = \|\hat{h}_+\|^2 + \|\hat{h}_-\|^2 .
\ee 
\end{thm}

\section{Symmetry and kernel structures of 2D FMT, FT, QFT, QFMT}

The QFMT of real signals analyzes symmetry. 
The following notation will be used\footnote{In this section we assume $g\neq \pm f$, but a similar study is possible for $g=\pm f$.}. 
The function $h_{ee}$ is \textit{even} with respect to (w.r.t.)
$r \rightarrow \frac{1}{r} 
\Longleftrightarrow 
\ln r \rightarrow -\ln r$,
i.e. w.r.t. the reflection at the unit circle,
and \textit{even} w.r.t. $\theta \rightarrow -\theta$,
i.e. w.r.t. reversing the sense of rotation
(reflection at the $\theta = 0$ line of polar coordinates
in the ($r,\theta$)-plane). Similarly we denote by
$h_{eo}$ even-odd symmetry, by 
$h_{oe}$ odd-even symmetry, and by 
$h_{oo}$ odd-odd symmetry.   

\begin{figure}
\begin{center}
  \includegraphics[height=7.0cm]{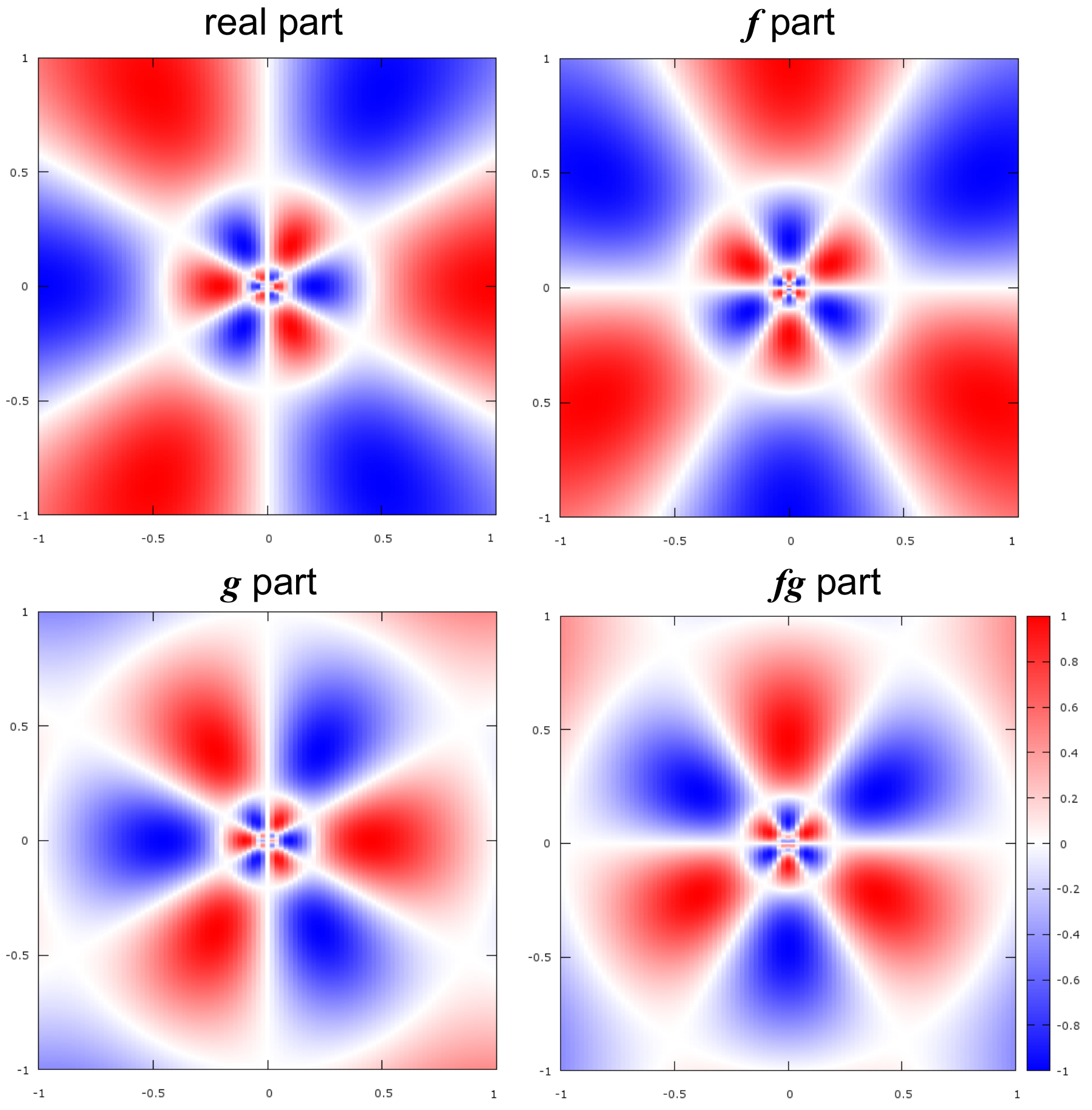}
  \includegraphics[height=7.0cm]{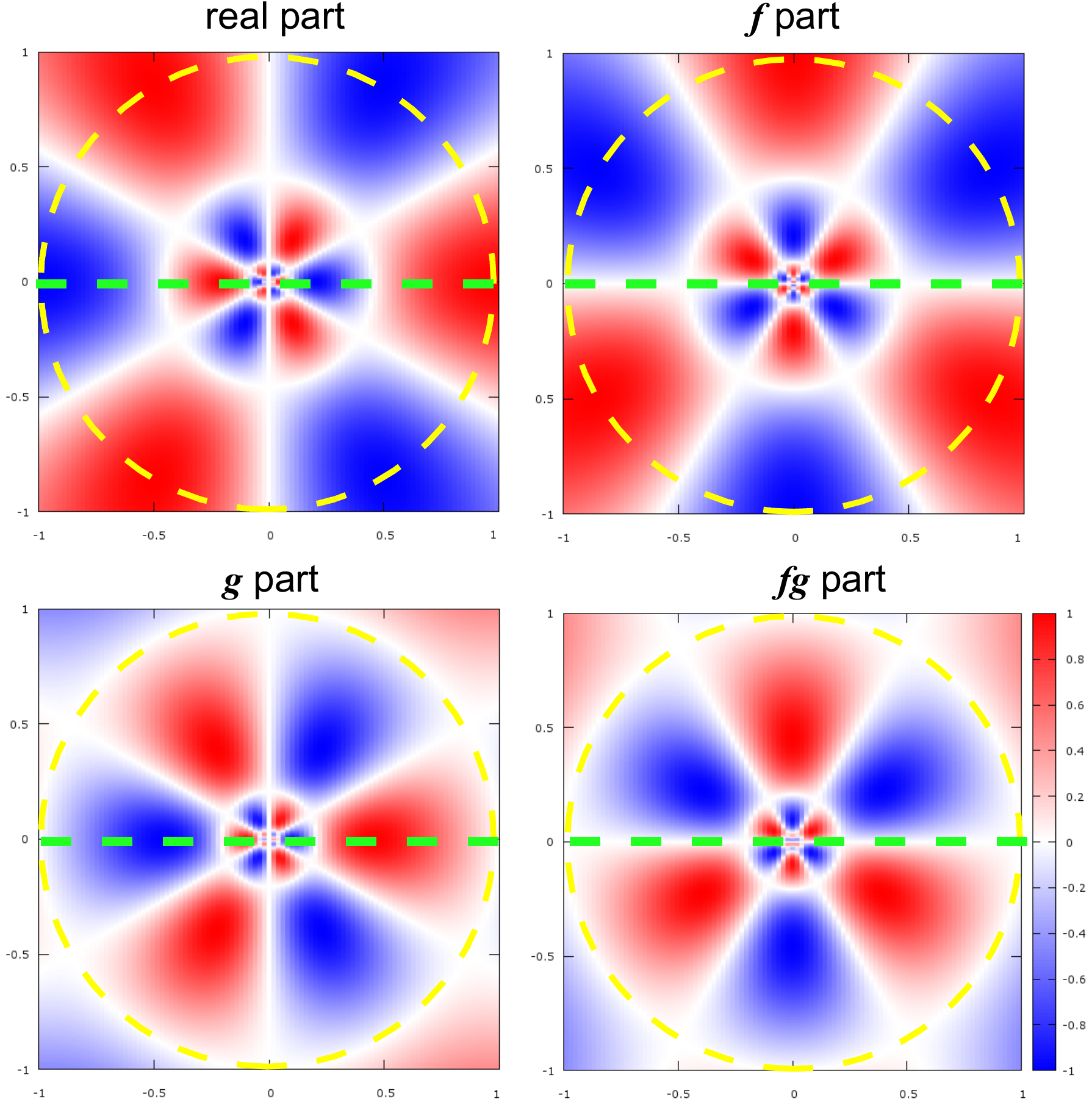}
  \caption{\textit{Left:} Four components of the QFMT kernel ($v=2, k=3$). 
  \textit{Right:} Symmetries of four components of the QFMT kernel ($v=2, k=3$). 
  \label{fg:4cQFMTsym}}  
\end{center}
\end{figure}

Let $h$ be a real valued function $\R^2 \rightarrow \R$. The QFMT of $h$ results in
\be
  \hat{h}(v,k) 
  = \underbrace{\hat{h}_{ee}(v,k)}_{\text{real part}}
    + \underbrace{\hat{h}_{eo}(v,k)}_{f\text{-part}}
    + \underbrace{\hat{h}_{oe}(v,k)}_{g\text{-part}}
    + \underbrace{\hat{h}_{oo}(v,k)}_{fg\text{-part}} \, .
\ee 
The QFMT of a real signal therefore automatically separates components
with different combinations of symmetry w.r.t. reflection at the unit circle and
reversal of the sense of rotation.
The four components of the QFMT kernel differ by radial and angular phase shifts, see the left side of Fig. \ref{fg:4cQFMTsym}.
The symmetries of {$r\rightarrow 1/r$} (reflection at \textit{yellow unit circle}), 
and {$\theta \rightarrow -\theta$} (reflection at \textit{green line}) can be clearly seen on the right of Fig. \ref{fg:4cQFMTsym}.

\begin{figure}
\begin{center}
  \includegraphics[height=5.5cm]{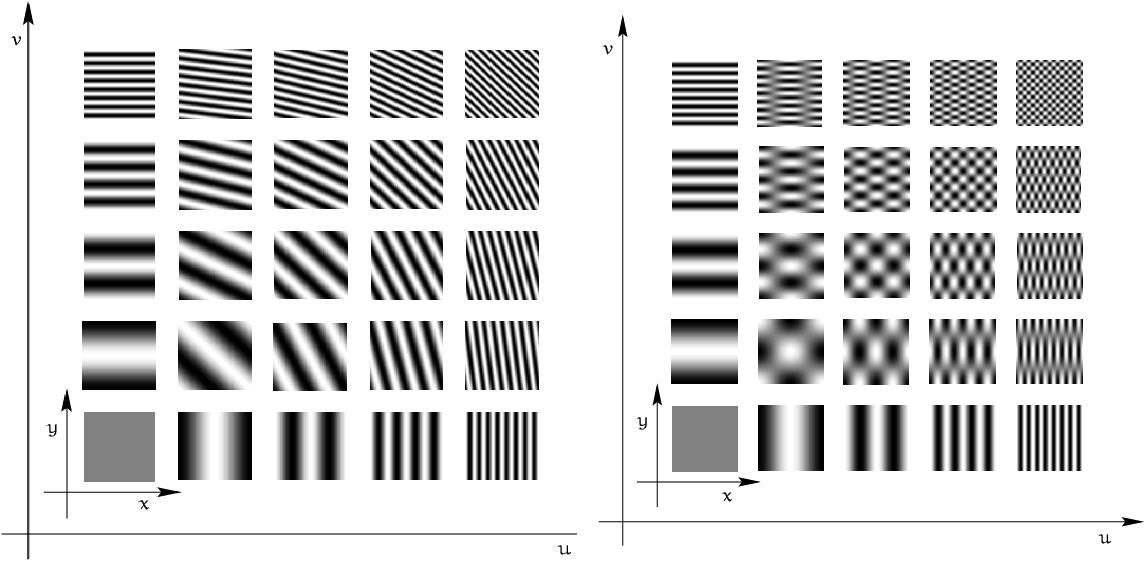}
  \caption{Left: 2D FT is intrinsically {1D}. Right: QFT is intrinsically {2D}. Source: \cite{TB:thesis}.
  \label{fg:kerFT+QFT}}
\end{center}
\end{figure}

\begin{figure}
\begin{center}
  \includegraphics[height=3.7cm]{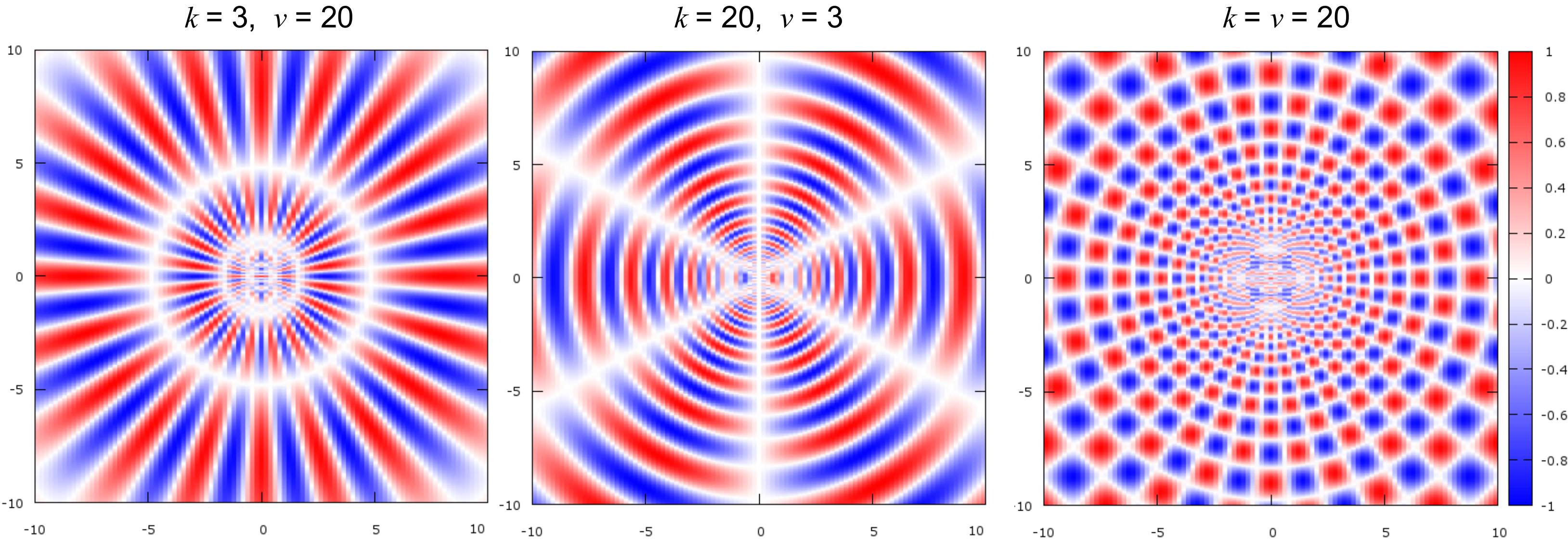}
  \caption{
    \textit{Left}: High angular resolution. 
    \textit{Center}: High radial resolution. 
    \textit{Right}: High radial and angular resolution. 
    \label{fg:QFMTvk}
}
\end{center}
\end{figure}

\begin{figure}
\begin{center}
  \includegraphics[height=3.5cm]{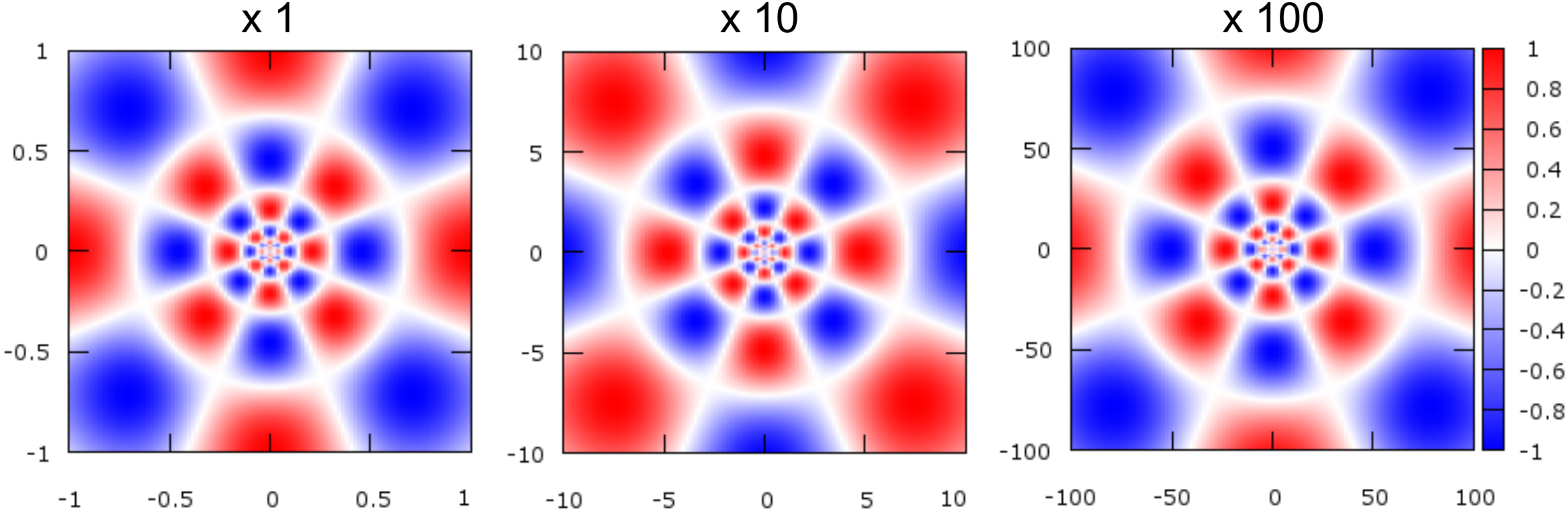}
  \caption{Illustration of QFMT scaling. \label{fg:QFMTscales}}
\end{center}
\end{figure}

Figure \ref{fg:QFMTvk} shows real the component of the QFMT kernel for various values of $v,k$, demonstrating various angular and radial resolutions. 
Figure \ref{fg:QFMTscales} shows the real component of the QFMT kernel for $v=k=4$ at three different scales. Similar patterns appear at all scales.
Figure \ref{fg:kerFT+QFT} shows the kernels of complex 2D Fourier transform (FT) $e^{-i(ux+vy)}$, $i \in \C$,
and the QFT $e^{-\sqi ux} e^{-\sqj vy}$, $\qi, \qj \in \HQ$, taken from 
~\cite{TB:thesis}, which treats applications to 2D gray scale images. Corresponding applications to 
color images can be found in \cite{ES:HFTColIm}. The 2D FT is intrinsically 1D, the QFT is intrinsically 2D, which makes it superior in disparity estimation and 2D texture segmentation, etc.  

\begin{figure}
\begin{center}
  \includegraphics[height=5.5cm]{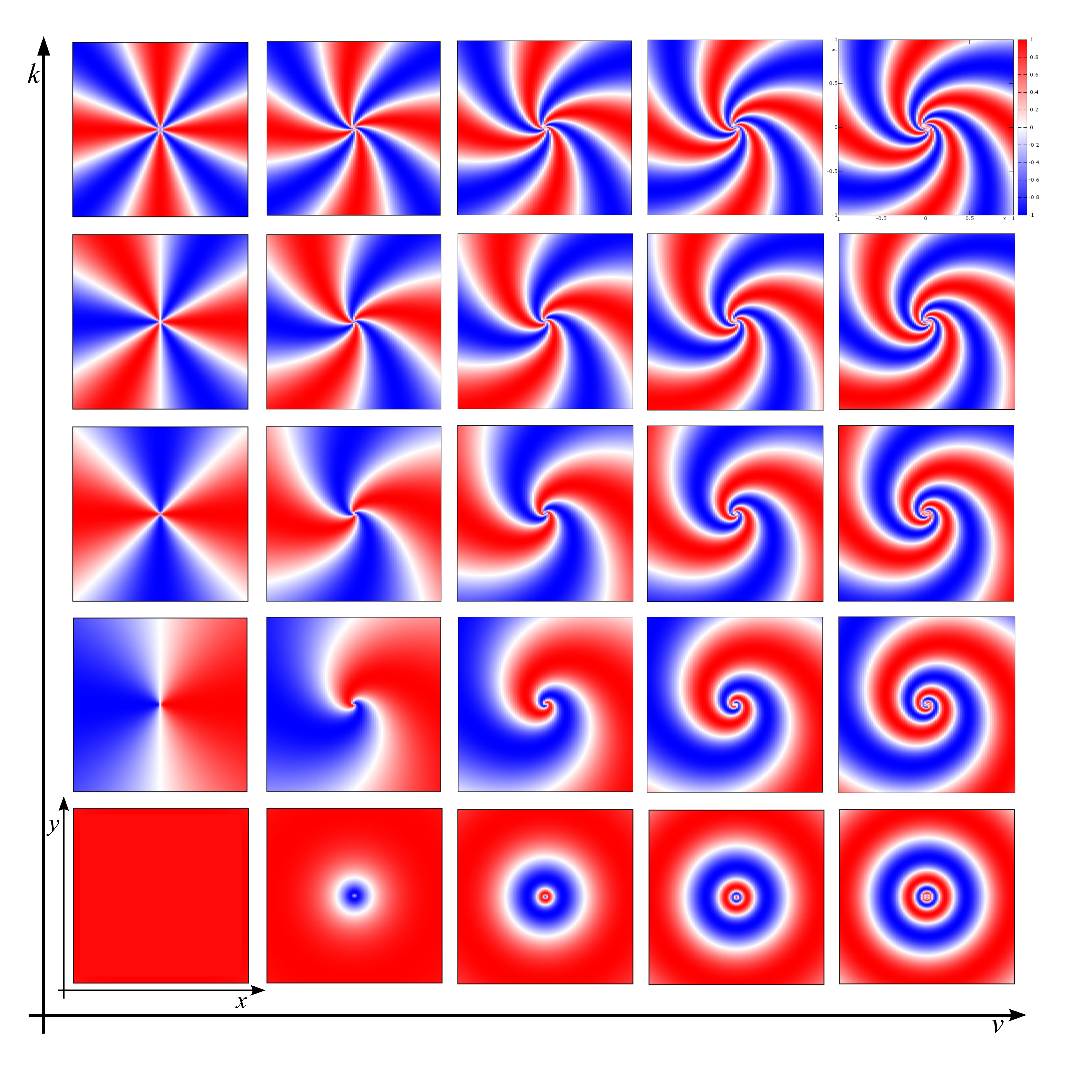}
  \includegraphics[height=5.5cm]{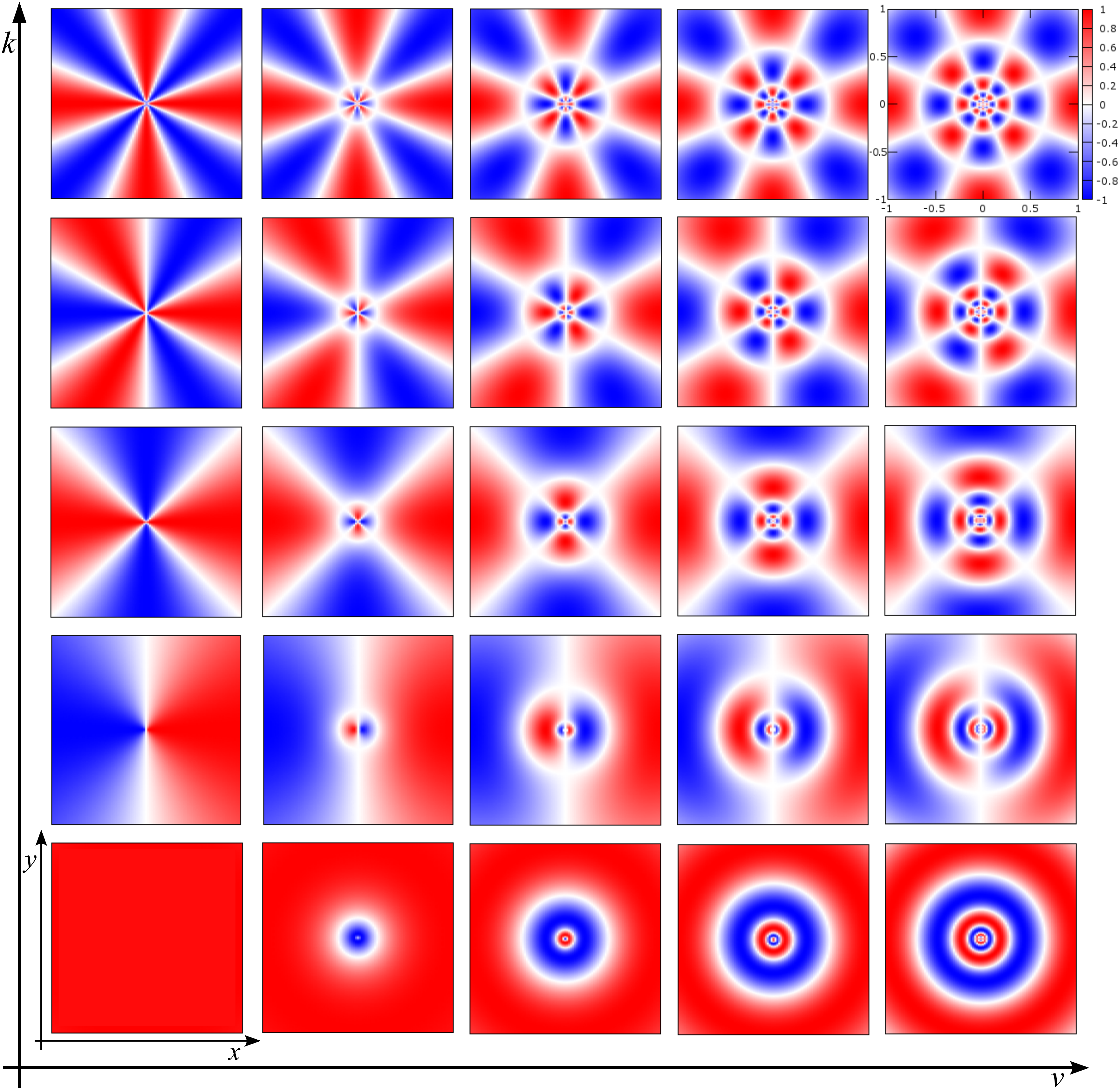}
  \caption{\textit{Left}: Kernel of FMT. 
  \textit{Right}: Kernel of QFMT. 
  $k,v \in \{0,1,2,3,4\}.$
  \label{fg:FMT+QFMT}}
\end{center}
\end{figure}

\begin{figure}
\begin{center}
  \includegraphics[height=5.5cm]{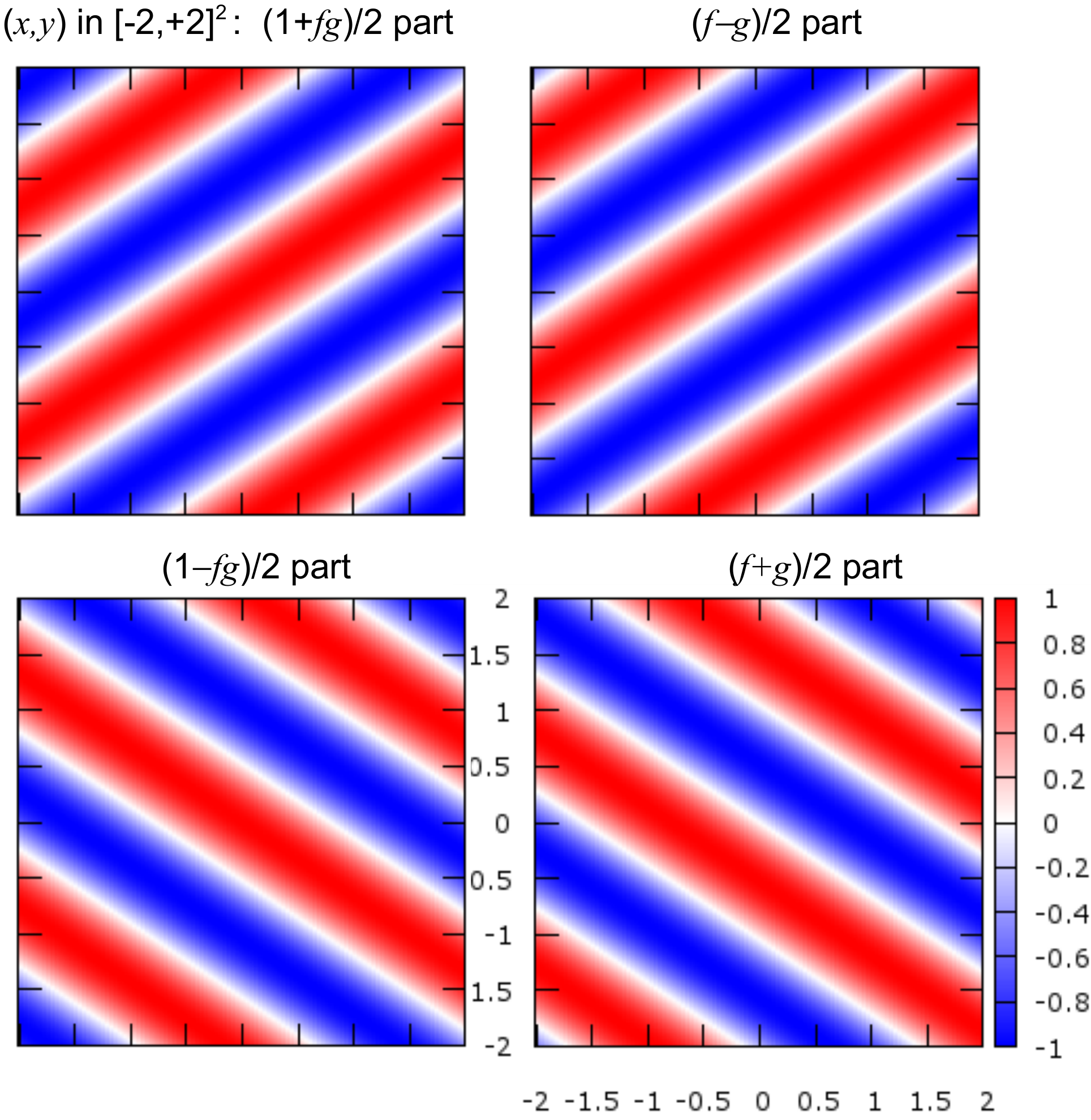}
  \includegraphics[height=5.5cm]{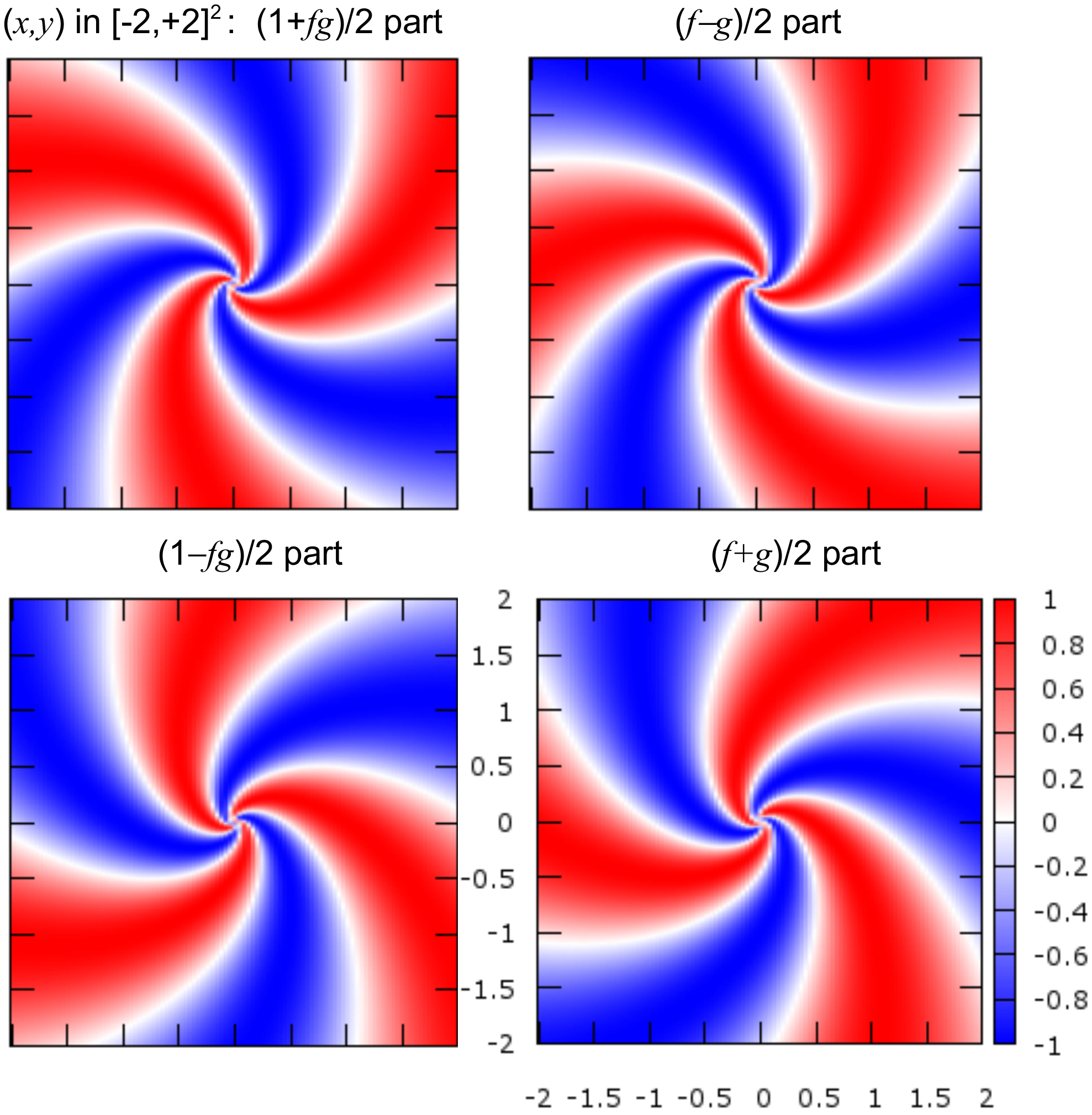}
  \caption{Left: QFT kernel, right: QFMT kernel.  
  \textit{Top row}:    $q_+$ parts: $1+fg$ and $f-g$ components. 
  \textit{Bottom row}: $q_-$ parts: $1-fg$ and $f+g$ components.
  \label{fg:QFT+QFMT}}
\end{center}
\end{figure}

Figure \ref{fg:FMT+QFMT} compares the kernels (real parts) of 2D complex FMT and the QFMT. Obviously the 2D QFMT can analyze genuine 2D textures better than the 2D complex FMT. 
Finally, Fig. \ref{fg:QFT+QFMT} compares the {kernels of the QFT} (left) and the {QFMT} (right). The scale invariant feature of the QFMT is obvious. Compared with the left side of Fig. \ref{fg:FMT+QFMT}, the QFMT is obviously the linear superposition of two quasi-complex FMTs with opposite winding sense, as shown in Theorem \ref{th:fpmtrafo}. 

\section{Conclusion}

The algebra of quaternions allows to construct a variety of 
{quaternionic Fourier-Mellin transformations} (QFMT), dependent on the choice of $f,g \in \HQ$, $f^2=g^2=-1$. Further variations would be to place both kernel factors initially at the left or right of the signal $h(r,\theta)$. The whole QFMT concept can easily be generalized to Clifford algebras $Cl(p,q)$, based on the general theory of square roots of $-1$ in $Cl(p,q)$. 

The modulus of the transform is scale and rotation invariant. Preceeded by 2D FT or by QFT, this allows {translation, scale and rotation invariant object description}. 
A diverse range of {applications} can therefore be imagined: Color object shape recognition, color image registration, application to evaluation of hypercomplex integrals, etc. 
  
Future research may be on extensions to Clifford algebras $Cl(p,q)$, to windowed and wavelet transforms, discretization, and numerical implementations.

\section*{Acknowledgments}

\begin{quote}
Zechariah after the birth of his son John:
And you, my child, ... will go on before the Lord to prepare the way for him,
to give his people the knowledge of salvation
   through the forgiveness of their sins,
because of the tender mercy of our God,
   by which the rising sun will come to us from heaven
to shine on those living in darkness
   and in the shadow of death,
to guide our feet into the path of peace. 
(Bible, Luke 1:76-79)
\end{quote}
 
I thank my wife, my children, my parents, as well as
B. Mawardi, G. Scheuermann, R. Bujack, G. Sommer, S. Sangwine, Joan Lasenby, 
H. Ishi, T. Sugawa, K. Matsuzaki. 
In memoriam: H. Shitaka (Hiroshima, 1945--2009), G. Major (1969--2011).

\section*{The Creative Peace License}
\label{pg:CPL}
\begin{enumerate}
\item
The study and application of this research is only permitted for peaceful,  non-offensive and non-criminal purposes. This permission includes passive devensive technologies, like missile defense shields. 
\item
Any form of study, use and application of research, which is licensed under the \textit{Creative Peace License}, for military purposes, with the explicit or implicit intent to create or contribute to military offensive technologies is strictly prohibited.
\item 
Individuals, groups, teams, public and private entities engaging in any form in the study, use and application of research licensed under the \textit{Creative Peace License}, thereby agree in a legally binding sense to strictly adhere to the terms of the \textit{Creative Peace License}. 
\end{enumerate}

\newpage

\end{document}